\documentclass[11pt]{article}

\title{Cayley sum graphs and eigenvalues of $(3,6)$-fullerenes}

\author{
   Matt DeVos
 \and
   Luis Goddyn\thanks{Supported by a Canada NSERC Discovery Grant}
 \and
   Bojan Mohar\thanks{Supported in part by ARRS Research Grant P1--0297, by
    an NSERC Discovery Grant, and by the CRC program.}~\thanks{On leave from:
  IMFM \& FMF, Department of Mathematics, University of Ljubljana, Ljubljana,
  Slovenia.}
 \and
   Robert \v{S}\'{a}mal\thanks{Supported by PIMS postdoctoral fellowship.}
   \thanks{On leave from Institute for Theoretical
     Computer Science (ITI), Charles University, Prague, Czech Republic.}
  \\
  \\
  {Department of Mathematics}\\
  {Simon Fraser University}\\
  {Burnaby, B.C. V5A 1S6} \\
  email: {\tt \{mdevos,goddyn,mohar,rsamal\}@sfu.ca}
}

\date{}

\usepackage{fullpage}
\usepackage{graphicx}
\usepackage{amsfonts}
\usepackage{amsmath}
\usepackage{amsthm}
\usepackage{hyperref}

\theoremstyle{plain}
\newtheorem{theorem}{Theorem}[section]
\newtheorem{lemma}[theorem]{Lemma}

\newtheorem{construction}[theorem]{Construction}

\theoremstyle{definition}
\newtheorem{example}[theorem]{Example}

\newcommand\CayS{\mathrm{CayS}}
\newcommand\Cay{\mathrm{Cay}}
\newcommand\Z{\mathbb{Z}}

\newcommand\R{\mathbb{R}}
\newcommand\E{\mathbb{E}}
\newcommand\T{\mathcal{T}}
\newcommand\HH{\mathcal{H}}
\newcommand\proj{\mathbf{p}}
\newcommand\Abar{\bar A}
\newcommand\Bbar{\bar B}
\newcommand\Cbar{\bar C}

\newcommand\calA{\mathcal {A}}
\newcommand\calU{\Lambda_{\bullet}}
\newcommand\calD{\Lambda_{\circ}}
\newcommand\calL{\Lambda}
\newcommand\ba{\mathbf {a}}
\newcommand\bb{\mathbf {b}}
\newcommand\bc{\mathbf {c}}
\newcommand\bu{\mathbf {u}}
\newcommand\bv{\mathbf {v}}

\newcommand\myvector[2]{\binom{#1}{#2}}

\newcommand\diag{\mathop{\mathrm{diag}}}
\let\eps\varepsilon

\begin{document}

\maketitle


\begin{abstract}
We determine the spectra of cubic plane graphs whose faces have sizes 3 and 6.
Such graphs, ``(3,6)-fullerenes'', have been studied by chemists 
who are interested in their energy spectra.
In particular we prove a conjecture of Fowler, which
asserts that all their eigenvalues come in pairs of the form
$\{\lambda,-\lambda\}$ except for the four eigenvalues $\{3,-1,-1,-1\}$.
We exhibit other families of graphs which are ``spectrally nearly bipartite''
in this sense.
Our proof utilizes a geometric representation to recognize the algebraic 
structure of these graphs, which turn out to be examples of Cayley sum graphs.
\end{abstract}



\section{Introduction}\label{sec:intro}

A \emph{$(3,6)$-fullerene} is a cubic plane graph whose faces have
sizes 3 and 6. These graphs have received recent attention from
chemists due to their similarity to ordinary fullerenes. (Such
graphs are sometimes called \emph{$(3,6)$-cages} in that community, but in
graph theory this term already has a different, well-established
meaning.) In 1995, Patrick Fowler (see \cite{FJS}) conjectured the
following result, which we prove here. Prior to this work, this
result had been established for several subfamilies of
$(3,6)$-fullerenes \cite{FJS,Ceul,SMZ}. Recall that
the \emph{spectrum} of a graph is the multiset of eigenvalues
of its adjacency matrix.

\begin{theorem}
\label{main} If $G$ is a $(3,6)$-fullerene, then the spectrum of $G$
has the form $\{3,-1,-1,-1\} \cup L \cup -L$ where $L$ is a multiset
of nonnegative real numbers, and $-L$ is the multiset of their
negatives.
\end{theorem}

\begin{figure}
  \centerline{\includegraphics[width=10cm]{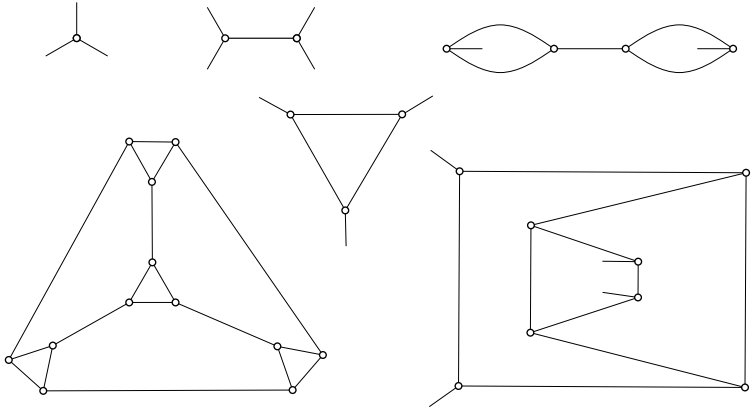}}
  \caption{Examples of some small $(0,3,6)$-fullerenes.}
  \label{fig:examples}
\end{figure}

In fact we prove an extended conjecture of Fowler \textit{et al.}\
\cite{FJS}. They propose that a generalized class of graphs called
\emph{$(0,3,6)$-fullerenes} also exhibit this ``spectrally nearly
bipartite'' behavior. A \emph{semiedge} of a graph is an edge with
one endpoint, but unlike a loop, a semiedge contributes just one to
both the valency of its endpoint and the corresponding diagonal
entry of the adjacency matrix.  In a plane embedding, a semiedge $s$
with endpoint $v$ is drawn as an arc with one end at $v$ which sits
in a face $f$, and $s$ contributes one to the length of $f$. 
Figure~\ref{fig:examples} displays some examples of small
$(0,3,6)$-fullerenes.

The outline of our proof is as follows. We show that every
$(0,3,6)$-fullerene can be represented as a quotient of a certain
lattice-like graph in the plane.  This geometric description allows
us to prove that these graphs are Cayley sum graphs.  Then we call
on a theorem which describes the spectral behavior of Cayley sum
graphs in terms of the characters of the group.

In fact, the geometric description of $(0,3,6)$-fullerenes which is
inherent in our proof is just a slight extension of a construction
for $(3,6)$-fullerenes which has been discovered by several authors
\cite{Ceul, FJS, Thur}, and follows easily from a deep theorem
on the intrinsic metric of polygonal surfaces by Alexandrov \cite{Alexandrov}.
In Section~\ref{sec:geom},
we give a proper exposition of this construction, and a proof that
it is universal.

With this construction in hand, it is possible to explicitly compute
the spectrum of $(0,3,6)$-fullerenes, and in Section
\ref{sec:compute} we detail precisely how this computation can be
made.  Finally, in Section \ref{sec:further}, we generalize to an
arbitrary dimension, and expose the connection between this type of
geometric graphs and Cayley sum graphs.

\section{Cayley sum graphs}\label{sec:cayley}

Let $\Gamma$ be a finite additive abelian group, and let $S
\subseteq \Gamma$.
We define the \emph{Cayley sum graph} $\CayS(\Gamma,S)$ to be the graph
$(V,E)$ with $V=\Gamma$, and $uv \in E$ if and only if $u+v\in S$.
If $S$ is a multiset, then $\CayS(\Gamma,S)$ contains multiple
edges, and if there exists $u \in \Gamma$ with $2u \in S$, then the
edge $uu$ is a semiedge.  This definition is a variation of the
well-studied \emph{Cayley graph} $\Cay(\Gamma,S)$, in which $uv$
forms an edge if and only if $u-v \in S$.


In contrast with Cayley graphs,
there are only a few appearances of Cayley sum graphs in the literature
(see \cite{GLS} and references therein).
For this reason we state some of their elementary properties. The
graph $G = \CayS(\Gamma,S)$ is $|S|$-regular.
While $G$ is not generally vertex-transitive, the map $x
\mapsto x+t$ is an isomorphism from $G$ to $\CayS(\Gamma,S+2t)$, for
every $t \in \Gamma$. Finally, the squared graph $G^{(2)}$, which
has an edge for each walk of length 2 in $G$, is the ordinary Cayley
graph $\Cay(\Gamma,S-S)$ where $S-S$ is the multiset $\{s_1-s_2 \mid
s_1,s_2\in S\}$.


The spectrum of a (finite abelian) Cayley graph~$\Cay(\Gamma,S)$
is easy to describe (see \cite[Ex.~11.8]{CPE} or \cite{L-sp}, where
the nonabelian case is dealt with). Every character $\chi$
of~$\Gamma$ is a (complex-valued) eigenvector corresponding to the
eigenvalue
\[
  \chi(S) := \sum_{s \in S} \chi(s) \,.
\]
We may assume
$\Gamma = \Z_{n_1} \times \dots \times \Z_{n_u}$,
where $|\Gamma| = \prod_i n_i$
and $\Z_k$ denotes the cyclic group of order~$k$.
To each $a = (a_1,\dots,a_u) \in \Gamma$
we associate the group character
\[
  \chi_a:(x_1,\dots,x_u) \mapsto \exp\left(2\pi i \sum_j \frac{a_j x_j}{n_j}\right).
\]
%
The characters for $a$ and $-a$ satisfy $\chi_{-a}(x)=
\overline{\chi_a(x)}$, so $\chi_a$ is a real-valued (indeed $\pm
1$-valued) eigenvector of $\Cay(\Gamma,S)$ if and only if $a$ is an
involutive group element. If $a$ is not involutive, then the real
and imaginary parts of $\chi_a$ provide real-valued eigenvectors for
the conjugate pair of eigenvalues $\chi_a(S), \chi_{-a}(S)$.

Cayley sum graphs exhibit a similar phenomenon.
Let $R = \{\chi_a \mid a+a=0\}$
be the real-valued characters of~$\Gamma$,
and let $C$ be a set containing exactly one character
from each conjugate pair
$\{\chi_a, \chi_{-a}\}$ (where $a\in \Gamma$ and $a+a\neq 0$).
So the set of characters of~$\Gamma$ is $R \cup \{ \chi,
\overline\chi \mid \chi \in C \}$. Versions of the following result
can be found in the literature~\cite{Chung, Alon}.


\begin{theorem} \label{eigenvalues}
Let $G = \CayS(\Gamma,S)$ be a Cayley sum graph on a finite abelian group $\Gamma$,
and let $R$, $C$ be as above.
The multiset of eigenvalues of~$G$ is
\[
   \{ \chi(S) \;:\; \chi \in R \}  \cup  \{ \pm |\chi(S)|\;:\;\chi \in C \} .
\]
The corresponding eigenvectors are $\chi$ (for $\chi \in R$),
and the real and the imaginary parts of\/ $\alpha \chi$
(for $\chi \in C$ with a suitable complex scalar $\alpha$
which depends only on $\chi(S)$).
\end{theorem}

\begin {proof}
Let $\chi$ be a character of~$\Gamma$ and $u \in \Gamma$ a vertex of $\CayS(\Gamma,S)$. Then
\[
  \sum_{v \in N(u)} \chi(v) = \sum_{s\in S} \chi(s-u) = \chi(S) \overline{\chi(u)} \,.
\]
This shows that every real-valued character is an eigenvector
corresponding to the eigenvalue~$\chi(S)$. If $\chi \in C$, then
$\chi$ is not an eigenvector. In this case we choose a complex
number $\alpha$ such that $|\alpha|=1$ and $\alpha^2 \chi(S) =
|\chi(S)|$ and we define $x(v) = \alpha \chi(v)$. It follows that
for every $u\in\Gamma$,
\[
  \sum_{v \in N(u)} x(v)
  = \alpha^2 \chi(S) \cdot \alpha^{-1} \overline{\chi(u)}
  = |\chi(S)| \cdot \overline{x(u)}.
\]
Consequently,
$\mathop{\mathrm{Re}} x$ and $\mathop{\mathrm{Im}} x$
are real eigenvectors corresponding
to eigenvalues $|\chi(S)|$ and~$-|\chi(S)|$, respectively.
Both of these vectors are nonzero, as they generate the same 2-dimensional
(complex) vector space as $\{\chi, \overline\chi\}$.
This, together with the orthogonality of characters,
implies that we have described the complete set of eigenvectors, and
thus the entire spectrum of~$\CayS(\Gamma,S)$.
\end {proof}


\section{$\mathbf {(0,3,6)}$-fullerenes as Cayley sum graphs}
\label{sec:caysum}

The goal of this section is to prove that $(0,3,6)$-fullerenes
are Cayley sum graphs, and to subsequently prove
Fowler's conjecture regarding their spectra.



The proof of Theorem~\ref{cayley} utilizes structural properties of
3-regular hexangulations of the torus. This class of graphs was
classified by Altshuler~\cite{A} and studied by many others
(e.g., Thomassen~\cite{Thom}). In a recent work of Alspach and Dean~\cite{ad},
it is shown that they are indeed Cayley graphs, and a description of the group is given.
Although the properties we require of these graphs are similar to those
found elsewhere, our approach is novel since it is inherently geometric.


A \emph{polygonal surface} $\HH$ is a connected 2-manifold without
boundary which is obtained from a collection of disjoint simple
polygons in $\mathbb{E}^2$ by identifying them along edges of equal
length. Thus we view $\HH$ both (combinatorially) as an embedded
graph with vertices, edges, and faces, and as a manifold with a
(local) metric inherited from~$\mathbb{E}^2$.

\begin{theorem}\label{cayley}
Every $(0,3,6)$-fullerene is isomorphic to a Cayley sum graph for an abelian
group which can be generated by two elements.
\end{theorem}

\begin{proof}
Let $G$ be a cubic $(0,3,6)$-fullerene with vertex set $V$.  Let
$G_2 = G \times K_2$ (the categorical graph product). Let
$(V_\bullet,V_\circ)$ be the corresponding bipartition of $V(G_2)$,
and for every $v \in V$, let $v_\bullet\in V_\bullet$ and $v_\circ\in
V_\circ$ be the vertices of $G_2$ which cover~$v$. Every semiedge
$vv \in E(G)$ lifts to the edge $v_\bullet v_\circ$
in~$G_2$. Each facial walk of $G$ bounding a face of size 6 lifts to
two closed walks of length 6 in~$G_2$, and each facial walk of $G$
bounding a face of size 3 lifts to a closed walk of length 6
in~$G_2$. Accordingly, we may extend $G_2$ to a polygonal surface
$\HH$ by treating all edges as having equal length and adding a
regular hexagon to each closed walk which is the pre-image of a
facial walk of $G$, with clockwise orientation as given by the
clockwise orientation of that face. Now, $\HH$ is an orientable
polygonal surface, all vertices have degree three, and all faces are
regular hexagons, so $\HH$ is a regular hexangulation of the flat
torus.
Let $\widetilde{\HH}$ be the
universal cover of $\HH$ and let $\proj : \widetilde{\HH}
\rightarrow \HH$ be the covering map.
Then $\widetilde{\HH}$ (with the metric inherited from $\HH$)
is the regular hexangulation of the Euclidean plane. We define
$\widetilde{V}_{\bullet} = \proj^{-1}(V_{\bullet})$,
$\widetilde{V}_{\circ} = \proj^{-1}(V_{\circ})$, and $\tilde{x}
= \proj^{-1}(x)$ for $x \in V_\bullet \cup V_\circ$.

Fix a vertex $u_{\bullet} \in V_{\bullet}$, and treat $\widetilde{\HH}$
as a regular hexangulation of $\mathbb{E}^2$ with $\proj((0,0)) =
u_{\bullet}$. This equips $\widetilde{\HH}$ with an (additive abelian)
group structure. The point set $\widetilde{V}_{\bullet}$ is a
geometric lattice. The point set $\tilde{u}_{\bullet}$ is a sublattice
of $\widetilde{V}_{\bullet}$.
Any fundamental parallelogram of $\tilde{u}_{\bullet}$
is a fundamental region of the cover~$\proj$.
We may identify $\HH$ with $\widetilde{\HH} / \tilde{u}_{\bullet}$, and
this equips $\HH$ with a group structure whose identity is
$u_{\bullet}$.

%

For every $y \in \HH$ ($y \in \widetilde{\HH}$) the map $x \mapsto x
+ y$ is an isometry of $\HH$ ($\widetilde{\HH}$, respectively). This
map may or may not preserve the combinatorial structure of~$\HH$
($\widetilde{\HH}$). 
An isometry $\mu : \HH \rightarrow \HH$
is \emph{respectful} if $\mu$ is an automorphism of
the embedded graph associated with $\HH$. An isometry $\tilde{\mu} :
\widetilde{\HH} \rightarrow \widetilde{\HH}$ is \emph{respectful} if
it is a lift of a respectful isometry of $\HH$.
%
%
Now, for every $y \in \widetilde{V}_{\bullet}$ the map $x \mapsto x+y$
is a respectful isometry of $\widetilde{\HH}$. Accordingly,
$V_{\bullet}$ is a subgroup of $\HH$ with identity $u_\bullet$, and for
every $y \in V_{\bullet}$ the map $x \mapsto x+y$ is a respectful
isometry of~$\HH$.

Let $\rho$ be the automorphism of the graph $G_2$ given by the rule
$\rho(v_{\circ}) = v_{\bullet}$ and $\rho(v_{\bullet}) = v_{\circ}$
for every $v \in V$.  Now, $\rho$ extends naturally to a respectful
isometry of $\HH$ which preserves the orientation of the hexagons,
but interchanges $V_{\circ}$ and~$V_{\bullet}$. We choose a
respectful isometry $\tilde{\rho}$ of $\mathbb E^2$ so that $\rho$
lifts to $\tilde{\rho}$.
Because $\tilde\rho$ preserves the orientation of $\mathbb E^2$, the
isometry $\tilde{\rho}$ is either a rotation or translation. Since
$\tilde \rho$ is respectful and maps $\widetilde{V}_\bullet$ to
$\widetilde{V}_\circ$, it easily follows that $\rho$ is either a
rotation by $\pi$ about the center of an edge or face, or $\rho$ is
a rotation by $\pi/3$ about the center of a face $F$.

%

We first consider the latter case. Here, all three vertices of
$\widetilde{V}_\bullet$ which are on the boundary of $F$, lie in the
same orbit of $\tilde{\rho}^2$. Since $\rho^2$ is the identity, all
three vertices cover the same vertex, say $v_\bullet$ in $\HH$. The
other three vertices of $F$ cover~$v_\circ$. In this case $G_2$ is
the theta-graph with vertex set $\{v_\bullet,v_\circ\}$. Here we
have $G \cong \CayS( \{0\}, \{0,0,0\} )$, the graph with one vertex
and three semiedges, and there is nothing left to prove.

We henceforth assume that $\tilde{\rho}$ is a rotation by $\pi$. Let
$x,y \in V_{\bullet}$ and choose $\tilde{x}, \tilde{y} \in
\widetilde{V}_{\bullet}$ which project (respectively) to $x$, $y$.
Then (using the fact that $\tilde{\rho}$ is a rotation by $\pi$) we
find that
\begin{align*}
\rho ( \rho(x) + y)
    &=  \proj( \tilde{\rho}( \tilde{\rho}(\tilde{x}) + \tilde{y}))  \\
    &=  \proj( \tilde{x} - \tilde{y} )  \\
    &=  x - y.
\end{align*}
In other words, for any fixed $y \in V_{\bullet}$, conjugating the map
on $\HH$ given by $x \mapsto x+y$, by $\rho$ yields the map $x
\mapsto x-y$.

We define a labeling $\ell : V_{\bullet} \cup V_{\circ} \to V_\bullet$
by the rule $\ell(v_{\bullet}) = \ell(v_{\circ}) = v_{\bullet}$. We
regard $\ell$ to be a labeling of $V(G_2)$ by elements of the
abelian group $V_\bullet$.
Let $v \in V$ and let $y \in V_{\bullet}$.  Then we have
\begin{align*}
\ell(v_{\bullet} + y) &=  \ell(v_{\bullet}) + y \\
\intertext{and} \ell(v_{\circ} + y)
    &=  \ell( \rho( v_{\circ} + y ))     \\
    &=  \ell( \rho(  \rho(v_{\bullet}) + y ))   \\
    &=  \ell( v_{\bullet} - y ) \\
    &=  \ell( v_{\circ} )  - y.
\end{align*}
That is, the group $V_{\bullet}$ acts on the labels of points in
$V_{\bullet}$ by addition and on the labels of points in $V_{\circ}$
by subtraction. Let $S$ be the multiset of labels of the three
vertices in $V_{\circ}$ which are adjacent to $u_{\bullet}$ (recall
that $u_\bullet$ is the group identity for $V_\bullet$). Then, for every
$v_\bullet \in V_{\bullet}$,
the three neighbors of $v_\bullet$ in $G_2$ have labels $S - v_\bullet$.
In particular, $v$ and $v'$ are adjacent vertices in $G$ if and only
if $\ell(v_{\bullet}) + \ell(v'_{\circ}) \in S$.  It follows
immediately from this that $G \cong \CayS(V_{\bullet},S)$. Since
$\widetilde{V}_{\bullet}$ can be generated by two elements, $V_{\bullet}
= \widetilde{V}_\bullet/\widetilde u_\bullet$ can also be generated by
two elements, and this completes the proof.
\end{proof}


We need only one quick observation before we resolve
Theorem~\ref{main} and the extended conjecture of Fowler et al.  If
$G$ is a cubic plane graph with $s$ semiedges, and $f_i$ faces of
size $i$ for every $i \ge 1$, then $3|V(G)| = 2|E(G)| - s =
\sum_{i\ge 1} i f_i $. Applying Euler's formula, we find that
$\sum_{i\ge 1}(6-i)f_i=12-3s$. In particular, every
$(0,3,6)$-fullerene satisfies
\begin{equation}\label{euler}
s+f_3 = 4.
\end{equation}




\begin{theorem} \label{full}
If $G$ is a $(0,3,6)$-fullerene with $s$ semiedges, then the
spectrum of $G$ may be partitioned as $\{M,L,-L\}$ where one of
the following holds:

  {\rm (a)} $s=0$ and $M = \{3,-1,-1,-1\}$,
  
  {\rm (b)} $s=2$ and $M = \{3,-1\}$,
  
  {\rm (c)} $s=3$ and $M = \{3\}$, or
  
  {\rm (d)} $s=4$ and $M = \{3,1\}$.
\end{theorem}

\begin{proof}
By the previous theorem, there is an abelian group $\Gamma$ which
can be generated by two elements so that $G \cong \CayS(\Gamma,S)$
for some $S \subseteq \Gamma$ with $|S| = 3$. By Theorem
\ref{eigenvalues}, we may partition the eigenvalues of $G$ into
multisets $M,L,-L$ where $M = \{ \chi(S) : \chi \in R \}$ and $R$ is
the set of $\pm 1$-valued characters of $\Gamma$. Every eigenvalue
in $M$ is the sum of three integers in $\{\pm 1\}$. The identity
character corresponds to $3\in M$. Since $G$ is not bipartite, we
have $-3 \notin M$, so every other element of $M$ is $\pm 1$. The
trace of the adjacency matrix is equal to~$s$, and is also equal to
the sum of the eigenvalues. Since $L$ and $-L$ sum to~$0$, we
conclude that $s=\sum M$.

We have $|R| \in \{1,2,4\}$ because $\Gamma$ has $2^k$ involutive
elements, for some $k\le 2$. If $|R|=1$, then $M = \{3\}$ and $s=3$
as in the statement. If $|R|=2$, then $s=\sum M = 3 \pm 1$, so we
have either the case $s=2$ or $s=4$ of the statement. Finally, we
assume $|R|=4$. By Equation~\eqref{euler} we have $s \le 4$, so
$\sum M \in \{0,2,4\}$. If $\sum M = 0$, then $s=0$ ($G$ is a
$(3,6)$-fullerene), and we have case~(a).
Finally, if $\sum M \in\{2,4\}$, then $M$ contains both a $1$ and
a~$-1$. By transferring these two entries from $M$ to the multisets
$L$ and $-L$, we find ourselves again in either the case $s=2$ or
the case $s=4$ of the statement. This completes the proof.
\end{proof}

We remark that there are infinitely many $(0,3,6)$-fullerenes with $s$ semiedges,
for each $s = 0$, $2$, $3$, $4$. As shown by Theorem~\ref{full}, there
are none with $s=1$, a fact that is non-trivial to prove from the first principles
(compare Theorem~2 (with $k=3$) in \cite[Sec.~13.4, p.~272]{Grunbaum}).

\section{An explicit construction}
\label{sec:geom}

It is known that all $(3,6)$-fullerenes
arise from the so-called \emph{grid construction}.
Roughly speaking, the grid construction expresses the dual plane graph,
which is a triangulation of the sphere,
as a quotient of the regular triangular grid.
The grid construction is also used by physicists \cite{Ceul,SMZ}
(sometimes without formal justification)
since it a convenient way to classify $(3,6)$-fullerenes
and compute their invariants.

We describe an extension of the grid construction
and show that it characterizes the $(0,3,6)$-fullerenes.
The construction makes clear how semiedges arise.
The group structure of $(0,3,6)$-fullerenes
is explicitly determined as a quotient of the $A_2$ lattice group.
With this, we can easily find the Cayley sum graph representation
via standard lattice computations,
and thereby determine the spectrum and the eigenvectors
of every $(0,3,6)$-fullerene.

In the following, $\T$ is the infinite triangular grid, whose
vertices (called \emph{gridpoints}) are points in the $A_2$ lattice.
The midpoint of any edge in $\T$ is called an \emph{edgepoint}.
The \emph{dual} $G^*$  of a plane graph $G$ with semiedges
is defined as an obvious extension of the dual of an ordinary graph;
every semiedge in $G$ which is incident with vertex $v$ and face $f$
corresponds to a semiedge in $G^*$ which is incident with the
dual vertex $f^*$ and the dual face $v^*$.

\begin{construction}
\label{constr}
The following procedure results in a
$(0,3,6)$-fullerene, $G$.
\begin{enumerate}
\item
Let $\triangle ABC$ be a triangle having no obtuse angle, and whose
vertices are gridpoints of $\T$. Let $\Abar$, $\Bbar$, $\Cbar$ be
the midpoints of the edges which are opposite to $A$, $B$, $C$
(respectively) in $\triangle ABC$.
\item
Optionally, translate $\triangle ABC$ so that $A$ coincides with an
edgepoint of $\T$.
\item From $\triangle ABC$, we fold an (isosceles) tetrahedron $Q = A\Abar
\Bbar \Cbar$ by identifying the boundary segment $\Abar B$ with
$\Abar C$, $\Bbar C$ with $\Bbar A$, and $\Cbar A$ with $\Cbar B$
(so $A$, $B$, and $C$ are identified into a single vertex in $Q$).
The portion of $\T$ lying within $\triangle ABC$ becomes a finite
graph $G^*$, possibly with semiedges, and drawn on the surface of
$Q$.
\item
Let $G$ be the dual of the plane graph $G^*$.
\end{enumerate}
\end{construction}

Every gridpoint within or on the boundary of $\triangle ABC$, except
$A$, $\Abar$, $\Bbar$, and $\Cbar$, has degree $6$ in $G^*$, and
corresponds to a hexagonal face of~$G$. After Step~2, each of $A$,
$\Abar$, $\Bbar$, $\Cbar$ is either a gridpoint or an edgepoint
of~$\T$. If $X \in \{A, \Abar, \Bbar, \Cbar\}$ is a gridpoint, then
$X$ becomes a vertex of degree $3$ in $G^*$, and corresponds to a
triangular face in $G$. If $X$ is an edgepoint, then $X$ becomes one
end of a semiedge in $G^*$, which corresponds to a semiedge in $G$.
It follows that Construction~\ref{constr} results in a
$(0,3,6)$-fullerene.

We remark that Construction~\ref{constr} works even if $\triangle ABC$ has
an obtuse angle (although it does not yield a geometric tetrahedron).
However, this does not give any new $(0,3,6)$-fullerenes, as the following
theorem shows. By forbidding obtuse triangles, we lose no generality
and gain canonicality. 

\begin {figure}
  \centerline{\includegraphics[width=15cm]{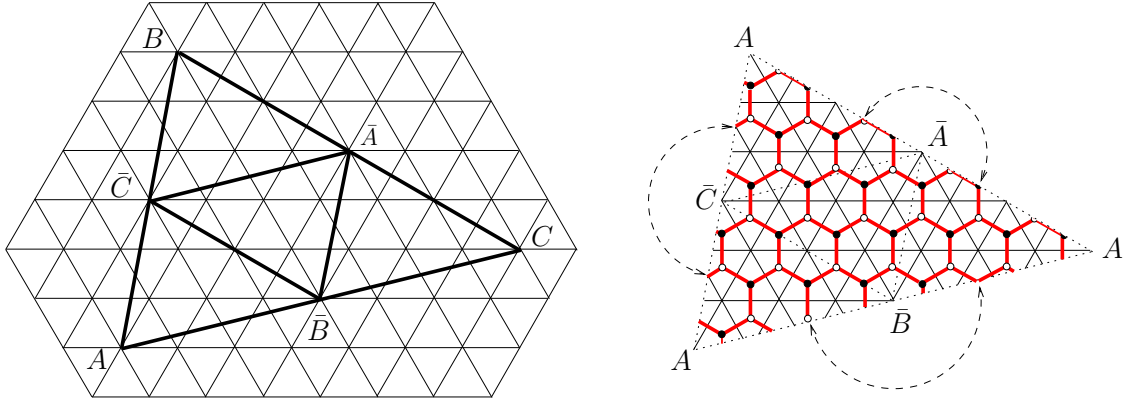}}
  \caption{An example of Construction~\ref{constr}.}
  \label{fig:grid}
\end {figure}

\begin{theorem} \label{thm:constr}
Every $(0,3,6)$-fullerene arises from Construction~\ref{constr}.
\end{theorem}
\begin {proof}
Let $G$ be a $(0,3,6)$-fullerene. If $G$ has just one vertex, then
$G$ arises from the construction when $\triangle ABC$ is a
triangular face of $\T$. We assume $G$ has at least two
vertices. The proof of Theorem \ref{cayley} shows that the direct
product $G_2 = G \times K_2$ is a bipartite hexangulation $\mathcal
H$
of the flat torus, where $\mathcal H$ is the image of a covering map
$\proj: \widetilde{\mathcal H}\to \mathcal H$ from a hexagonal
tessellation of the plane.

We further recall that there is an isometry $\rho$ of $\mathcal H$
which is respectful of $G_2$ and interchanges its partite sets
$V_\bullet$ and $V_\circ$. This isometry lifts to an isometry of
$\widetilde{\mathcal H}$ which is a rotation $\tilde\rho$ by  $\pi$
about a point, say $A\in \widetilde{\mathcal H}$, which is either
the center of a hexagonal face, or the midpoint of an edge of
$\widetilde{\mathcal H}$.
(More precisely $\tilde\rho:x \mapsto 2A-x$ is the central symmetry
through $A$.) The kernel of $\proj$ is a geometric lattice $\Lambda$
in $\widetilde{\mathcal H}$, and rotation by $\pi$ about any point
in in the scaled lattice
$\frac12 \Lambda$ also projects to $\rho$.
Let $B$, $C$ be points in $\widetilde{\mathcal H}$ such that the
vectors $AB$, $AC$ form lattice basis for $\Lambda$. By possibly
translating $C$ by a (unique) integer multiple of $AB$, we can
assume that $\triangle ABC$ has no obtuse angles. This lattice basis
defines a fundamental parallelogram $ABDC$ where $AD = AB+AC$.
Scaling the parallelogram by $\frac12$ results in a fundamental
parallelogram for
$\frac12 \Lambda$ whose vertices we may label $A \Cbar \Abar \Bbar$
as in Construction~\ref{constr}.

Now each vertex $v$ of $G$ lifts to a unique pair of vertices
$v_\bullet, v_\circ$ in the (half-open) parallelogram $ABDC$, where
$v_\bullet, v_\circ$ are centrally symmetric about $\Abar$. By
representing $v$ with the vertex in $\{v_\bullet, v_\circ\}$ which
lies in $\triangle ABC$, and identifying boundary segments of
$\triangle ABC$ as in Step~3 of the construction, we obtain an
isomorphic copy of $G$. Finally, Construction~\ref{constr} is stated
in terms of the triangular grid $\T$, which is the plane dual of
$\widetilde{\mathcal H}$.
\end {proof}

We remark that Construction~\ref{constr} in fact produces a
$(0,3,6)$-fullerene $G$ rooted at a triangle or semiedge labeled
with $A$. Two triangles drawn in $\T$ result in isomorphic
pairs $(G,A)$ if and only if the triangles are congruent. Therefore
the map $\triangle ABC \mapsto G$ is at most 4-to-1 up to symmetries
of $\T$.


\section{Computing the spectrum}
\label{sec:compute}

In this section, we use Construction~\ref{constr} to compute the
group and spectrum of any particular $(0,3,6)$-fullerene~$G$.

The faces of $\T$ consist of \emph{up-triangles} ($\Delta$)
and \emph{down-triangles} ($\nabla$).
We regard $\calU$ to be an
$A_2$-lattice generated by unit-length vectors $\ba, \bb$ with
$\angle\ba\bb=\pi/3$. With $A$ being the gridpoint selected in
Step~1 of Construction~\ref{constr},
we shall assume that the origin of $\calU$
is (the center of) the up-triangle $u_{\bullet} := \triangle A (A+\ba) (A+\bb)$.
Note that $\calU$ is a translation of
the gridpoints of $\T$ and corresponds to $\widetilde V_\bullet$ in
the proof of Theorem \ref{cayley}.
We denote by $\calL$ the sublattice of $\calU$ generated by vectors~$\overrightarrow{AB}$
and~$\overrightarrow{AC}$.
(A translation of $\calL$ is used in the proof
of Theorem~\ref{thm:constr}.)
%
In Step~2, we translate $\triangle ABC$ by a vector
\begin{equation}\label{p}
\bc := \frac {p_1}2 \ba + \frac {p_2}{2}\bb
\end{equation}
for integers $p_1$, $p_2$.
We may assume without loss of generality that $p_1, p_2 \in \{0,1\}$,
so, after Step~2,
the point $A$ is
either a vertex or an edgepoint on the boundary of $u_\bullet$.
Let $p$, $q$, $r$, $s$ be integers satisfying
\begin{equation}\label{pqrs}
AB = p\ba+q\bb, \qquad AC = r\ba+s\bb.
\end{equation}
Let $\Abar$, $\Bbar$, $\Cbar$ and $T$ be as in the construction of~$G$.

To express $G$ as a Cayley sum graph
we label the faces of $\T$ with elements of the finite abelian group
presented as $\Gamma = \langle \,\alpha, \beta \,\mid\,
  p\alpha + q\beta = 0, \,
  r\alpha + s\beta = 0 \,\rangle$.
We define $f:\calU \to \Gamma$ by
\begin{equation}\label{f}
f(i\ba+j\bb)=i\alpha+j\beta,
\end{equation}
and extend $f$ to the down-triangles in such a way that triangles which are
centrally symmetric with respect to $A$ receive the same value of $f$.
The kernel of $f$ is (a translation of) the lattice $\calL$ generated
by $AB$ and $AC$.
We observe the following properties:
\begin {itemize}
  \item $f$ assigns the same value to triangles that are
    identified during the `folding' stage of the construction.
    This is because the triangles that are identified are symmetric with
    respect to each of $\Cbar$, $\Bbar$, and $\Abar$;
    each of these symmetries is a composition of
    the negation map $\calU \mapsto -\calU$
    and a translation by an element of
    $\calL = \ker f$.
  \item $f$ is a bijection from $V(G)$ to $\Gamma$.
     By construction, the up-triangles within the fundamental region
     $ABDC$ correspond to elements of~$\Gamma$.
     The down-triangles within the triangle $ABC$
     correspond to up-triangles within~$DCB$.
  \item If $u_1$ and $u_2$ are two up-triangles, then
      $f(u_2) = f(u_1) + f(u_2 - u_1)$.
      If $d_1$ and $d_2$ are two down-triangles then
      $f(d_2) = f(d_1) - f(d_2 - d_1)$.
\end {itemize}
Now let $u$ be any up-triangle and $d_1$, $d_2$, $d_3$ its neighbors.
We define the \emph{sum-set} $S = \{f(u) + f(d_i) \mid i=1,2,3\}$.  From
the above-mentioned properties of~$f$ it follows that $S$ does not
depend on the choice of~$u$.
The symmetry around~$A$ shows that we get the same sum-set
if we consider neighbors of a down-triangle to define~$S$.
It follows that $G \cong \CayS( \Gamma, S)$.

We can explicitly compute $\Gamma$ and $S$
by applying standard lattice computations.
We recall that the \emph{Smith normal form} of a nonsingular
integer matrix $M$ is the unique matrix
$\diag(\delta_1,\delta_2,\dots,\delta_k) = UMV$
where $U$ and $V$ are unimodular and
the product $\delta_1\delta_2\cdots\delta_i$ is the g.c.d.\ of
the order $i$ subdeterminants of $M$, $1 \le i \le  k$
(see, e.g. \cite[Section 4.4]{Sch-TLIP}).

\begin {lemma}   \label{cagesfla}
Let $G$ be a $(0,3,6)$-fullerene obtained from Construction~\ref{constr},
and let $\bc$, $p$, $q$, $r$, $s$ be as in \eqref{p} and \eqref{pqrs}.
Let $\diag(m,n) = UMV$ be the Smith normal form of the matrix
$M = \begin{pmatrix}
p & r \\
q & s \\
\end{pmatrix}$.
Let $\bu$, $\bv$ denote the columns of~$U$.
Then $G = \CayS(\Gamma, S)$ where $\Gamma = \Z_m \times \Z_n$ and
\[
S = \{\,(p_1-1)\bu +  p_2    \bv,\;
         p_1   \bu + (p_2-1) \bv,\;
        (p_1-1)\bu + (p_2-1) \bv\,\}.
\]
Here we interpret each column vector $\myvector{x_1}{x_2} \in S$
to be the group element $(x_1 \!\!\mod m, \,x_2 \!\!\mod n)\in\Gamma$.
\end {lemma}
\begin{proof}
The columns of the matrix $B := (\ba,\bb)$ form a lattice basis for $\calU$
whereas those of $BM$ generate the sublattice $\calL$.
Since $U$ and $V$ are unimodular,
the columns of $B' :=BU^{-1}$ also generate $\calU$.
Accordingly, $\calL$ is generated by the columns of $BMV = B'\,\diag(m,n)$.
It follows that $\Gamma = \calU / \calL \cong \Z_m \times \Z_n$.
If we index the up-triangles with respect to the basis $B'$,
then the mapping $f: B'\myvector{i'}{j'} \mapsto (i' \bmod m, j' \bmod n)$
is the one defined in \eqref{f}.
Changing the basis to $B = B'U$,
we find $f(i \ba + j \bb) = i \bu + j \bv$.

After Step~1 of the construction,
the three down-triangles which are neighbours of $u_\bullet$ reflect through $A$
to the up-triangles at $-\ba$, $-\bb$ and $-\ba-\bb$.
When $A$ is translated by $\bc$ in Step~2,
the three up-triangles are accordingly translated by $2\bc = p_1\ba + p_2\bb$.
Therefore
\[
S=\{f((p_1-1)\ba+p_2\bb),\, f(p_1\ba+(p_2-1)\bb),\, f((p_1-1)\ba+(p_2-1)\bb)\}
\]
as claimed.
\end{proof}

We present a sample computation illustrating the determination
of the group and spectrum.

\begin{example}
The example of Figure~\ref{fig:grid}
corresponds to $(p_1,p_2)=(0,0)$ and $(p,q,r,s)=(6,2,-2,6)$.
All six integers are even, so the resulting graph $G$
has no semiedges.
We compute the Smith normal form to be
\[
UMV
=\begin{pmatrix}
 0 &  1 \\
-1 & -7
\end{pmatrix}
\begin{pmatrix}
 6 & -2 \\
 2 & 6 \\
\end{pmatrix}
\begin{pmatrix}
-2 & -3 \\
 1 & 1 \\
\end{pmatrix}
=
\begin{pmatrix}
 2 & 0  \\
 0 & 20\\
\end{pmatrix}.
\]
Hence $\Gamma = \Z_2\times \Z_{20}$.
Furthermore, the generating set is
\[
S=\{ -\bu+0\bv, 0\bu-\bv, -\bu-\bv \}
 = \{(0,1),(-1,7),(-1,8)\}.
\]

This implies $G$ has eigenvalues $3, -1, -1, -1$, and
\[
\{ \pm | \eps^b + (-1)^a \eps^{7b} + (-1)^a \eps^{8b}|
\,:\,
    0 \le a \le 1, \,1 \le b \le 9 \} \,,
\]
where $\eps = e^{2\pi i/20}$.

If we were to translate $\triangle ABC$ by $(\frac12 \ba, 0)$,
then we get a $(0,3,6)$-fullerene $G'$ with four semiedges.
Here we have $(p_1,p_2)=(1,0)$, which has the effect of translating the
generating set by $\bu$.  That is,
\[
G' = \CayS(\:\Z_2\times \Z_{20},
\{(0,0),(1,6),(1,7)\}\,),\;\text{ and}
\]
\[
  \mathop{\mathrm{spec}}(G') = \{3,1,1,-1\} \cup
  \{ \pm | 1 + (-1)^a \eps^{6b} + (-1)^a \eps^{7b}|
\,:\,
    0 \le a \le 1,\, 1 \le b \le 9 \} \,.
\]
The four semiedges are incident with the vertices
$(0,0), (1,0), (0,10), (1,10)\in \Gamma$.
\end{example}

\section{The geometry of Cayley sum graphs}
\label{sec:further}

In Section~\ref{sec:geom} we saw how
the geometric description of $(0,3,6)$-fullerenes
in terms of the $A_2$ lattice
implies that they are Cayley sum graphs.
Therefore their eigenvectors are easy to calculate,
and their spectra are ``nearly bipartite.''
Here we describe the precise circumstances
under which Cayley sum graphs arise from
geometric lattices in this manner.
In fact we will see that \emph{every}
Cayley sum graph arises as a quotient
of two cosets of a geometric lattice.
We then exhibit some
families of Cayley sum graphs which
have a recognizable crystallographic local
structure.

First we review, in greater generality, the conditions under which a
graph $G$ is a Cayley sum graph. Let $G_2$ be the cover $G \times
K_2$ with bipartition $(V_{\bullet},V_{\circ})$. Note that $G_2$ has
a natural automorphism, $\rho$, -- we call it the
\emph{inversion map} -- which transposes the two vertices within
each fibre. By following the proof of Theorem~\ref{main}, we find
that $G$ is a Cayley sum graph on an abelian group $\Gamma$ if and
only if $\Gamma$ acts regularly on each of $V_\bullet$ and $V_\circ$
as a group of $G_2$-automorphisms, and this action satisfies
\begin{equation}\label{conj-cond}
  \rho^{-1}g\rho = -g \text{, \,for each $g \in \Gamma$.}
\end{equation}


Our construction proceeds with a sequence of graphs
\[
\widetilde G \mapsto \widetilde G_2 \mapsto G_2 \mapsto G.
\]
Let $\calU \subset \E^d$ be a geometric lattice, and let $\widetilde
G = \CayS(\calU,S)$ be a Cayley sum graph with edges drawn as
straight line segments. Each generator $s\in S$ corresponds to a set
of edges of $\widetilde G$ whose midpoints are concurrent at the
point $\frac12 s$. Let $\calD$ be any nontrivial coset of $\calU$,
and let $A\in\R^d$ be such that $\calD = 2A +\calU$. Let $\tilde
\rho: x \mapsto 2A-x$ be the inversion map through $A$. As above, we
construct $\widetilde G_2 = \widetilde G \times K_2$ with partite
sets $(\widetilde V_\bullet, \widetilde V_\circ) = (\calU, \calD)$,
where the fibres of $\widetilde G_2$ are the orbits of~$\tilde
\rho$.

The graph $\widetilde G_2$ is drawn in Euclidean $d$-space $\E^d$
with straight line segments for edges. Let $\E^d/ \tilde \rho$
denote the quotient space (an orbifold) whose points are the $\tilde
\rho$-orbits $\{x, \tilde \rho(x)\}$, $x \in \E^d$. Geometrically
speaking, $\E^d/ \tilde \rho$ is a cone with its apex $A$ having the
solid angle of a halfspace. By mapping each point in $\E^d$ to its
$\tilde \rho$-orbit, we may view $\widetilde G \cong \widetilde
G_2/\tilde \rho$ as being naturally embedded in $\E^d/ \tilde \rho$.
Every edge of $\widetilde G_2$ whose midpoint is $A$ folds to a
semiedge of $\widetilde G$. In the case of $(0,3,6)$-fullerenes,
$\widetilde{G}_2$ is the plane hexagonal grid, and $\widetilde{G}$
is a grid drawn on a cone where every face is a hexagon except at
$A$, where $A$ is either the midpoint of a triangular face, or the
end of a semiedge.

\begin{figure}
  \centerline{\includegraphics[width=15.5cm]{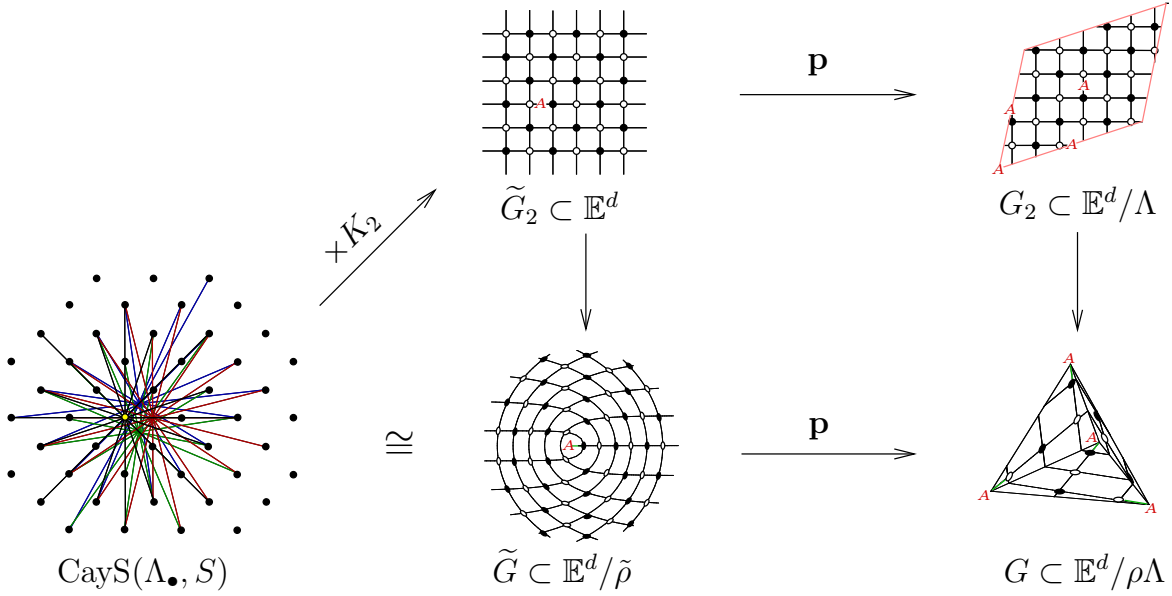}}
  \caption{Constructing finite Cayley sum graph from a lattice. 
  Illustrated with the $D_2$-lattice, resulting in a $28$-vertex
  Cayley sum graph which is also a $4$-regular quadrangulation of the tetrahedron.}
  \label{fig:commute}
\end{figure}

Now let $\calL$ be any sublattice of $\calU$, and let $\proj$ be the
be the natural projection from $\E^d$ to the $d$-torus $\E^d/\calL$.
Then $G_2 := \proj(\widetilde G_2)$ is a finite bipartite graph with
partite sets $(V_\bullet,V_\circ) := (\proj(\calU),\proj(\calD))$,
and which is embedded in $\E^d/\calL$. Then $\tilde \rho$ projects
to $\rho$, a symmetry of order $2$ in the $d$-torus. Evidently
$\rho$ is
an inversion map for $G_2$ satisfying~\eqref{conj-cond} with $\Gamma
= \calU/\calL$. Therefore $G \cong G_2/\rho$ is a finite Cayley sum
graph embedded in the orbifold $\E^d/\rho\calL$. Let $\calA \subset
\E^d/\calL$ be the fixed points of $\rho$. Then $\calA$ consists of
exactly $2^d$ points having the form $\proj(A + \frac12 \calL)$ and
$\rho$ acts on $\E^d/\calL$ as an inversion through any point in
$\calA$. As an orbifold, $\E^d/\rho\calL$ is orientable if and only
if $d$ is even. To visualize $\E^d/\rho\calL$, it is convenient to
select a fundamental region for $\E^d/\calL$ whose $2^d$ extreme
points belong to $A+\calL$. Let $T$ be the part of the region which
lies on the positive side of a hyperplane $H$, which contains the region's
centroid.
All points in $\calA$ lie on the boundary of $T$ so we
obtain $\E^d/\rho\calL$ by an appropriate gluing of the boundary of
$T$. The graph $G$ is embedded in $T$ with each vertex
$\{x,\rho(x)\}$ represented by the unique point in $\{x,\rho(x)\}
\cap T$. For example, $\E^2/\rho\calL$ is an isosceles tetrahedron,
whose four extreme points comprise~$\calA$. The grid construction of
$(0,3,6)$-fullerenes corresponds to selecting $H$ to be a diagonal
of a fundamental parallelogram. The Cayley sum graph $G$ has one
semiedge for every point of $\calA$ which lies on an edge of $G_2$.
Figure \ref{fig:commute} summarizes
the commuting projections
and the four embedded graphs.
%

Since every finite abelian group
is the quotient of two geometric lattices,
it follows that every finite Cayley sum graph $G$
arises from a
quadruple $(\calU, S, A, \calL)$
as described above.
By employing a linear transformation we can even
assume that $\calU = \Z^d$.
We do not make this assumption here, since that would
obfuscate the following examples.
When the sum set $S$
is a set of lattice points which are close to $2A$,
then each edge of $\widetilde G_2$ is a short line segment,
and $\widetilde G_2$ is often a recognizable
bipartite crystallographic configuration.
After selecting $\calL$ and applying the above construction,
we obtain a finite Cayley sum graph embedded in $T$
with a local geometry that reflects
the crystallographic structure of $\widetilde G_2$.
We present with some examples.
\begin{itemize}
\item
For $d=1$, if $\widetilde G_2$ is the two-way infinite path, then
$\widetilde G$ is the infinite ray with a semiedge at its origin,
and $G_2$ is an even cycle drawn as a regular polygon. The inversion
$\rho$ identifies points reflected in a line which bisects opposite
edges of the polygon, so $G$ is a finite path with a semiedge at
each end. The spectrum of $G$ takes the form $M \cup L \cup -L$
where $M=\{2\}$ or $M=\{2,0\}$.

\item \emph{(Grid-like examples)}
If $\calU = D_d$, the lattice of integer points of even weight, and
$\calD = \calL_{\bullet}+(1,0,0,\dots)$, then $\calU \cup \calD =
\Z^d$, and we may take $\widetilde G_2$ to be the standard cartesian
grid. If $A = (\frac12,0,0,\dots)$, then applying the construction
with any sublattice $\calL$ of $\calU$ leads to a Cayley sum graph
$G$ having exactly $2^d$ semiedges.

If $d=2$, then $G$ is a $4$-regular quadrangulation of an isosceles
tetrahedron, with a semiedge at each tetrahedral vertex.
Such an graph illustrates Figure~\ref{fig:commute}.  The set of
unmatched eigenvalues of $G$ is either $M=\{4\}$ or $M=\{4,0\}$.
Indeed, \emph{every} 4-regular quadrangulation of a sphere can be
expressed in this way. To see this fact,
we need only adapt the proof of Theorem~\ref{cayley}.

When $d$ is odd, we may take
$A = (\frac12,\frac12,\frac12,\dots)$.
Since $A$ is not on an edge of the cubic grid,
this results in a grid-like Cayley sum graph $G$
having fewer than $2^d$ semiedges.  Indeed
$G$ has no semiedges at all
if $\calL$ is a sublattice of $2\calU$.

\item \emph{(Diamond-like examples)}
Again we take $\calU$ to be the $D_d$-lattice, but put $\calD =
\calL_{\bullet} + (\frac12,\frac12,\frac12,\dots)$. The set $\calU
\cup \calD$ is commonly called the \emph{generalized diamond
packing}, and is denoted by $D_d^+$ (see \cite[p.~119]{Sloan}). The
\emph{diamond grid} is the graph $\widetilde{G}_2$ in which each
point in $\calL_\bullet$ is joined to the $2^{d-1}$ nearest points in
$\calD$. Putting $A =  (\frac14,\frac14,\frac14,\dots)$ results in a
Cayley sum graph having at least $2^{d-1}$ semiedges.
A more attractive option is to put $A =
(\frac54,\frac14,\frac14,\dots)$, which lies on no edge of
$\widetilde{G}_2$. Provided that $\calL$ is a sublattice of
$2\calU$, this results in a Cayley sum graph having no semiedges.
When $d=3$, this construction gives a class of Cayley sum graphs
having the local structure of diamond crystal. Such graphs
satisfy $M=\{4,-2,-2\}$ or $M=\{4,0,-2,-2\}$. Another attractive
class is based on $D_8^+$, otherwise known as the $E_8$ lattice.

\item
The 24-dimensional Leech lattice $\Lambda_{24}$
arises as the union of two cosets of
a lattice $h\Lambda_{24}$ which is obtained from
the binary Golay code (see \cite[p.\ 124]{Sloan}).
This yields a particularly attractive class of crystallographic
Cayley sum graphs of high dimension.
\end{itemize}

We have constructed infinite families of Cayley sum graphs whose
spectra have the form $M \cup L \cup -L$, where $M$ is a fixed finite multiset.
By taking the categorical products with a fixed graph $H$,
one obtains other ``spectrally nearly bipartite'' families of graphs.
It would be interesting to find other natural examples of this phenomenon.

\bibliographystyle{rs-amsplain}
\bibliography{CayleySum}

\end{document}